\documentclass[11pt]{amsart}
\usepackage{amsmath, amsthm, amscd, amssymb,hyperref}

\newtheorem{thm}{Theorem}[section]
\newtheorem*{thm*}{Theorem}
\newtheorem{cor}[thm]{Corollary}
\newtheorem*{cor*}{Corollary}
\newtheorem{lem}[thm]{Lemma}
\newtheorem{prop}[thm]{Proposition}

\newtheorem*{con*}{Conjecture}
\newtheorem*{prob*}{Problem}

\theoremstyle{definition}
\newtheorem{defn}[thm]{Definition}
\newtheorem*{ques}{Question}

\theoremstyle{remark}

\numberwithin{equation}{section}

\newcommand{\A}{\mathcal{A}}
\newcommand{\B}{\mathcal{B}}
\newcommand{\C}{\mathcal{C}}

\newcommand{\bbH}{\mathbb{H}}
\newcommand{\bbC}{\mathbb{C}}
\newcommand{\bbZ}{\mathbb{Z}}

\newcommand{\bbR}{\mathbb{R}}

\newcommand{\ii}{{\sqrt{-1}}}

\newcommand{\conf}{{\rm Conf}}
\newcommand{\gr}{{\rm gr}}

\DeclareMathOperator{\ad}{ad} \DeclareMathOperator{\Der}{Der}
\DeclareMathOperator{\Image}{Image} \DeclareMathOperator{\Aut}{Aut}
\DeclareMathOperator{\GL}{GL} \DeclareMathOperator{\IA}{IA}

\renewcommand{\b}[1]{\mathbf{#1}}

\begin{document}
\title[Centralizers of Lie Algebras]
{Centralizers of Lie Algebras
Associated to the Descending Central Series of Certain Poly-Free
Groups}

\author[D.~C.~Cohen]{Daniel C.~Cohen$^{*}$}
\address{Department of Mathematics, Louisiana State University,
Baton Rouge, LA 70803} \email{\href{mailto:cohen@math.lsu.edu}{cohen@math.lsu.edu}}
\urladdr{\href{http://www.math.lsu.edu/~cohen/}{www.math.lsu.edu/\~{}cohen}}
\thanks{$^{*}$Partially supported by National Security Agency grant
H98230-05-1-0055}

\author[F.~R.~Cohen]{F.~R.~Cohen$^{\dag}$}
\address{Department of Mathematics,
University of Rochester, Rochester, NY 14225}
\email{\href{mailto:cohf@math.rochester.edu}{cohf@math.rochester.edu}}
\thanks{$^{\dag}$Partially supported by the NSF}

\author[Stratos Prassidis]{Stratos Prassidis$^{\ddag}$}
\address{Department of Mathematics
Canisius College, Buffalo, NY 14208, U.S.A.}
\email{\href{mailto:prasside@canisius.edu}{prasside@canisius.edu}}
\thanks{$^{\ddag}$Partially supported by Canisius College Summer
Grant}



\begin{abstract}
Poly-free groups are constructed as iterated semidirect products of
free groups. The class of poly-free groups includes the classical
pure braid groups, fundamental groups of fiber-type hyperplane
arrangements, and certain subgroups of the automorphism groups of
free groups. The purpose of this article is to compute centralizers
of certain natural Lie subalgebras of the Lie algebra obtained from
the descending central series of poly-free groups $\Gamma$
including some of the geometrically interesting classes of groups
mentioned above. The main results here extend the result in
\cite{copr} for such groups. These results imply 
that a homomorphism $f: \Gamma \to {G}$ is
faithful, essentially, if it is faithful when restricted to the
level of Lie algebras obtained from the descending central series
for the product $F_T \times Z$, where $F_T$ is the ``top'' free group
in the semidirect products of free groups and $Z$ is the center of
$\Gamma$. The arguments use a mixture of homological, and Lie
algebraic methods applied to certain choices of extensions. The
limitations of these methods are illustrated using the ``poison
groups" of Formanek and Procesi \cite{fopr}, poly-free groups whose Lie
algebras do not have certain properties considered here.
\end{abstract}

\maketitle

\section{Introduction}

A group $\Gamma$ is {\it poly-free} if there is a sequence
\[
\{1\} = {\Gamma}_0 \le {\Gamma}_1 \le \dots \le {\Gamma}_n =
{\Gamma}
\]
such that ${\Gamma}_{i+1}$ is a normal subgroup of ${\Gamma}_i$ and
the quotient ${\Gamma}_i/{\Gamma}_{i+1}$ is isomorphic to a free
group. Examples include the Artin pure braid group, the fundamental
group of the configuration space of $n$ ordered points in $\bbC$, as
well as fundamental groups of certain orbit configuration spaces
\cite{mono,ckx}. There are also natural subgroups of the
automorphism group of a free group which are poly-free. One such
example considered here is the ``upper triangular" version of
McCool's  subgroup of basis-conjugating automorphisms \cite{copa}.
Homological properties of such groups have been investigated in
 \cite{cs} and \cite{JMM}. Geometric properties of poly-free
groups were investigated in \cite{afr} and \cite{fr}.

The main objective of this paper is to determine the centralizers of
certain Lie subalgebras of the Lie algebras which arise from the
classical descending central series of the families of discrete groups
given above. These Lie algebras have been used in contexts
ranging from Vassiliev invariants of pure braids \cite{kohno1,
kohno2} to structures of simplicial groups and the loop space of the
$2$-sphere \cite{BCWW,cowu1}. Other applications of these Lie algebras
arise in the study of homotopy groups of higher dimensional knot spaces, 
see \cite{sinha}.

For nonempty subsets $U$ and $V$ of a (discrete) group $G$, let
$[U,V]$ denote the subgroup of $G$ generated by all commutators
$[u,v]=u v u^{-1} v^{-1}$ of elements $u \in U$ and $v \in V$.  Let
$G_n$ be the $n$-th stage of the descending central series of $G$,
defined inductively by $G_1=G$, and $G_{n+1}=[G_n,G]$ for $n\ge 1$,
and let $\gr_n(G)=G_n/G_{n+1}$ be the $n$-th associated quotient.
Let $\gr_*(G)=\bigoplus_{n\ge 1} \gr_n(G)$. There is a bilinear
homomorphism
\[
[-,-]: {\gr}_p(G)\otimes_{\mathbb Z}{\gr}_q(G)\longrightarrow {\gr}_{p+q}(G)
\]
induced by the commutator map
\[
c:G \times G \longrightarrow G
\]
(which is itself not in general a homomorphism).

The construction of ${\gr}_*(G)$ is a functor from the category of
discrete groups to the category of Lie algebras over $\mathbb Z$.
This construction has potential applications to questions concerning
representations.  Specifically, if $G$ is {\it residually nilpotent}
and a group homomorphism $f$ out of $G$ induces a Lie algebra
monomorphism on the corresponding Lie algebras, then $f$ is itself a
monomorphism (\cite{cowu1}).  Based on this property, conditions
which insure that a representation of $G$ is faithful are recorded
in Corollary \ref{cor:embeddings} below.
It should be noted, however, that these methods 
have not yet succeeded in producing faithful representations of the
poly-free groups considered here.

\subsection{Summary of the Results}
Let $B$ be a subset of the Lie algebra $A$. Then {\it the
centralizer of $B$ in $A$ } is defined by
\[
{\C}_A(B) = \{a\in A:\; [a, b] = 0, \; \text{for all}\; b\in B\}.
\]
Abbreviate ${\C}_A(A)$ by ${\C}(A)$. The main result of this paper
gives centralizers of Lie subalgebras of $\gr_*(G)$ for the following groups $G$:

\begin{enumerate}
  \item[\textbf{I.}] The ``upper triangular McCool groups'',
subgroups $P\varSigma_n^+$ of the
     automorphism group of the free group $F_n=F[x_1,\dots,x_n]$ generated
by
     automorphisms $\beta_{i,j}$, $1\le i<j\le n$, defined by
\[
\beta_{i,j}(x_k) =
\begin{cases}
x_k,
  &   \text{if $k \neq i$, } \\
x_i^{-1}x_j^{}x_i^{},
  & \text{if $k =  j$. } \\
\end{cases}
\]
       \item[\textbf{II.}] Fundamental groups $P_G(n)$ of orbit
configuration spaces
     associated to surface groups $G$ acting freely, and properly
discontinuously
     on the upper $1/2$-plane $\mathbb H $. In this case,
${\gr}_*(P_G(n))$ is generated by $B_{i,j}^{\sigma}$ with $\sigma
\in G$, $1 \leq i < j \leq n$.
        \item[\textbf{III.}] Fundamental groups $P(r,n)$ of orbit configuration
spaces
     associated to finite cyclic groups, of order $r$,
acting freely on $\mathbb C^*$, by rotation. The generators of the
Lie algebra ${\gr}_*(P(r,n))$ are given by $B_{i,j}^{(p)}$ with
     $1 \le p \le r$, $1 \leq i < j \leq n$, and
     $Z_k$ for $1 \leq k \leq n$.
\end{enumerate}

The groups ${\Gamma}_n$ given by $P\varSigma_n^+$, $P_G(n)$, or
$P(r,n)$ all have the property that there is an extension
\begin{equation} \label{eq:split exact}
\begin{CD}
1 @>{}>>  F_{\alpha(n)}@>{j}>> {\Gamma}_n @>{p}>> {\Gamma}_{n-1}
@>{}>>  1
\end{CD}
\end{equation}
which satisfies the following properties.
\begin{subequations} \label{eq:conditions}
\begin{align}
&\text{$\Gamma_0$ is a free group;\hfill} \label{one}\\
&\text{$F_{\alpha(n)}$ is a free group on a countable set
$\alpha(n)$ of cardinality at least two;\hfill}
\label{two}\\
&\text{the map $p: {\Gamma}_n  \to \Gamma_{n-1}$ is a split
epimorphism;
and\hfill} \label{three}\\
&\text{the action of ${\Gamma}_{n-1}$ on $H_*(F_{{\alpha}(n)})$ is
trivial.} \label{four}
\end{align}
\end{subequations}
By work in \cite{kohno3,falk-rand,xi}, for such a group $\Gamma_n$,
there is a split extension of Lie algebras
\begin{equation} \label{eq:lie split exact}
\begin{CD}
0 @>{}>> {\gr}_*(F_{\alpha(n)}) @>{{\gr}_*(j)}>> {\gr}_*({\Gamma}_n)
@>{{\gr}_*(p)}>> {\gr}_*({\Gamma}_{n-1}) @>{}>>  0.
\end{CD}
\end{equation}
In particular the map $j: F_{\alpha(n)} \to {\Gamma}_n$ induces an
injection on the level of Lie algebras
\[
{\gr}_*(j): L[V_{\alpha(n)}]= {\gr}_*(F_{\alpha(n)}) \longrightarrow {\gr}_*(\Gamma_n),
\]
where $L[V_{\alpha(n)}]$ denotes the free Lie algebra generated by
the set ${\alpha}(n)$. Groups ${\Gamma}_n$ which admit this
structure are the main focus of this paper.

\begin{thm}\label{thm-main}
Let ${\Gamma}_n$ be one of the groups described above. The
centralizer of the free Lie algebra $L[V_{\alpha(n)}]$ in
${\gr}_*({\Gamma}_n)$ is given as follows:
\begin{enumerate}
\item[\textbf{I.}]  If ${\Gamma}_n =
P\varSigma_n^+$ is the upper triangular McCool group, then
\[
{\C}_{{\gr}_*({\Gamma}_n)}(L[V_{\alpha(n)}]) = L[B_n],
\] where $B_{i,j}$ is the image of $\beta_{i,j}$ in
$\gr_*(P\varSigma_n^+)$, and ${B}_n =
\sum_{j=2}^n{B}_{1,j}$.
\item[\textbf{II.}]  If ${\Gamma}_n = P_G(n)$ is the fundamental group of
a surface group orbit configuration space, then
\[
{\C}_{{\gr}_*({\Gamma}_n)}(L[V_{\alpha(n)}]) = 0.
\]
\item[\textbf{III.}]  If ${\Gamma}_n = P(r,n)$ is the fundamental group of a
cyclic group orbit configuration space, then
\[
{\C}_{{\gr}_*({\Gamma}_n)}(L[V_{\alpha(n)}]) = L[{\Delta}(r,n)],
\]
where
\[
{\Delta}(r,n) = \sum_{k=1}^nZ_k + \sum_{p=1}^r\sum_{1\le i < j \le
n}B_{i,j}^{(p)}.
\]
\end{enumerate}
\end{thm}

For a poly-free group $\Gamma_n$ satisfying conditions \eqref{eq:conditions}, 
the center of the Lie algebra $\gr_*(\Gamma_n)$ is determined in 
section \ref{sec:centers}.  Namely, it is shown in Proposition \ref{prop:cyclic centers} 
that the center of $\gr_*(\Gamma_n)$ is either infinite cyclic or trivial.  Furthermore, 
if the free group $\Gamma_0$ is of 
rank greater than one (as is the case for the groups $F_{\alpha(k)}$ by \eqref{two}), 
then both the center of $\gr_*(\Gamma_n)$ and the center of $\Gamma_n$ itself 
are trivial.  This property need not be satisfied for an arbitrary poly-free
group, as illustrated by the direct product $F_1 \times F_1$, where
$F_1$ is free of rank one.

Recall the classical adjoint representation
\[
\ad:L \longrightarrow \Der(L)
\]
of a graded Lie algebra $L$, where $\Der(L)$ denotes the graded Lie
algebra of graded derivations of $L$. The map $\ad$ is defined by
the equation $\ad(X)(Y) = [X,Y]$ for $X$ and $Y$ in $L$. 
The center of $L$ is precisely the kernel of $
\ad:L \to \Der(L)$. If $I$ is a Lie ideal of $L$, then the natural
restriction map is denoted by
\[
\ad|_I:L \longrightarrow \Der(I).
\]

Let ${\Gamma}_n$ be one of the groups considered above. Regard
${\gr}_*({\Gamma}_n)$ as a graded Lie algebra by the convention that
${\gr}_q({\Gamma}_n)$ has degree $2q$ (as the axioms for a graded
Lie algebra are not satisfied without this convention).  In these
cases, $L[V_{{\alpha}(n)}]$ is a Lie ideal. Consider the two adjoint
representations:
\[
\ad: {\gr}_*({\Gamma}_n) \longrightarrow \Der({\gr}_*(\Gamma_n)) \quad\text{and}
\quad
 {\ad}|_{L[V_{{\alpha}(n)}]}:
{\gr}_*({\Gamma}_n) \longrightarrow \Der(L[V_{{\alpha}(n)}]).
 \]

\begin{cor}\label{cor:kernel of adjoint representation}
Let ${\Gamma}_n$ be one of the above groups.
\begin{enumerate}
\item[\textbf{I.}] If ${\Gamma}_n = P\varSigma_n^+$, then
$\ker(\ad) = \ker({\ad}|_{L[V_{{\alpha}(n)}]})$ is the cyclic group
generated by ${B}_n$ in ${\gr}_1(P\varSigma_n^+)$. Thus there is
a short exact sequence of Lie algebras
\[
\begin{CD}
0 @>{}>> L[{B}_n] @>{}>>
{\gr}_*(P\varSigma_n^+)@>{\ad|_{L[V_n]}}>> \Image(\ad|_{L[V_n]})
@>{}>> 0.
\end{CD}
\]
\item[\textbf{II.}]  If ${\Gamma}_n = P_G(n)$, then
$\ker(\ad) = \ker({\ad}|_{L[V_{{\alpha}(n)}]})$ is trivial. Thus
there is an isomorphism
\[
\begin{CD}
{\gr}_*(P_G(n))@>{\ad|_{L[V_n]}}>> \Image(\ad|_{L[V_n]}).
\end{CD}
\]
\item[\textbf{III.}]  If ${\Gamma}_n = P(r,n)$, then
$\ker(\ad) = \ker({\ad}|_{L[V_{{\alpha}(n)}]})$ is the cyclic group
generated by $\Delta(r,n)$ in ${\gr}_1(P(r,n))$. Thus there is a
short exact sequence of Lie algebras
\[
\begin{CD}
0 @>{}>> L[\Delta(r,n)] @>{}>> {\gr}_*(P(r,n))@>{\ad|_{L[V_n]}}>>
\Image(\ad|_{L[V_n]}) @>{}>> 0.
\end{CD}
\]
\end{enumerate}
\end{cor}

Recent work on the Isomorphism Conjecture \cite{fajo:is, fajo:ri-ge}
has renewed interest in detecting monomorphisms on discrete groups,
especially when the target space is a finite dimensional linear
group and the image is a discrete subgroup.  The linearity of the
classical Artin braid groups, and hence the Artin pure braid groups,
was established by Bigelow \cite{bigelow} and Krammer
\cite{krammer}.  Subsequently, Digne~\cite{digne} showed that 
Artin groups of crystallographic type are linear.  
This paper may be viewed as an attempt to develop
methods for detecting faithful finite dimensional representations of
poly-free groups, as illustrated by the following consequence of
Theorem \ref{thm-main} and the results of~\cite{copr}.

\begin{cor}\label{cor:embeddings}
Let ${\Gamma}_n$ be one of the poly-free groups $P\varSigma_n^+$,
$P_G(n)$, or $P(r,n)$, and $f: {\Gamma}_n \to {G}$ a 
homomorphism of groups.  If the morphism of Lie algebras
\[
{\gr}_*(f)|_{L[V_{\alpha(n)}]}: L[V_{\alpha(n)}] \longrightarrow  {\gr}_*({G})
\]
is a monomorphism and
\begin{enumerate}
\item[\textbf{I.}] ${\gr}_*(f)|_{L[{B}_n]}: L[{B}_n]\to
{\gr}_*({G})$ is a monomorphism when $\Gamma_n=P\varSigma_n^+$;
\item[\textbf{II.}] no further conditions when $\Gamma_n=P_G(n)$;
\item[\textbf{III.}] ${\gr}_*(f)|_{L[\Delta(r,n)]}: L[\Delta(r,n)]\to
{\gr}_*({G})$ is a monomorphism when $\Gamma_n=P(r,n)$,
\end{enumerate}
then $f$ is a monomorphism. In addition, the following two
statements are equivalent:
\begin{enumerate}
    \item[(i)] The map $f: {\Gamma}_n \to {G}$ is one-to-one.
    \item[(ii)] The maps of Lie algebras
\[
{\gr}_*(f)|_{L[V_{\alpha(n)}]}: L[V_{\alpha(n)}] \longrightarrow {\gr}_*(f({\Gamma_n}))
\]
is a monomorphism and
\begin{enumerate}
\item[\textbf{I.}] ${\gr}_*(f)|_{L[{B}_n]}: L[{B}_n]\to\
{\gr}_*(f(P\varSigma_n^+))$ is a monomorphism when
$\Gamma_n=P\varSigma_n^+$;
\item[\textbf{II.}] no further conditions when $\Gamma_n=P_G(n)$;

\item[\textbf{III.}] ${\gr}_*(f)|_{L[\Delta(r,n)]}: L[\Delta(r,n)]\to
{\gr}_*(f(P(r,n))$ is a monomorphism when $\Gamma_n=P(r,n)$,
\end{enumerate}
where $f(-)$ denotes the image of $f$.
\end{enumerate}
\end{cor}

The Artin pure braid groups, and the groups $P(r,n)$, may be realized 
as fundamental groups of complements of fiber-type hyperplane arrangements.  
For any such arrangement, the fundamental group of the complement is poly-free, and satisfies the conditions \eqref{eq:conditions}, see Falk and Randell \cite{falk-rand}.  Thus, there are corresponding split, short exact sequences of descending central series Lie algebras as
in \eqref{eq:lie split exact}.  Furthermore, using the linearity of the pure braid group \cite{bigelow,krammer} and topological properties of fiber-type arrangements \cite{mono,ccx}, one can show 
that the fundamental group of the complement of any fiber-type arrangement is linear, see Theorem 
\ref{thm:ft linear}.  Consequently, many poly-free groups which fit into exact sequences of the form
\eqref{eq:split exact} admit faithful finite dimensional linear representations.

There is, however, a dichotomy as follows. The ``poison group''
$H$ of Formanek and Procesi admits no faithful, finite dimensional
linear representation \cite{fopr}. While this group is poly-free,
and fits into a split exact sequence
\[
\begin{CD}
1 @>{}>>  F_3 @>{j}>>  H @>{p}>> F_2 @>{}>>  1,
\end{CD}
\]
where $F_k$ is a rank $k$ free group, information is lost upon
passage to the descending central series Lie algebra. Specifically,
there are homomorphisms out of $H$ which have non-trivial kernels,
but induce monomorphisms on the level of descending central series
Lie algebras.
The natural map
\[
p\times \alpha:H \longrightarrow F_2 \times H_1(H),
\]
given by the projection $p:H\to F_2$ and the abelianization map
$\alpha:H \to H_1(H)$, is one such homomorphism. The map
$p{\times}{\alpha}$ has non-trivial kernel but induces a monomorphism on
the Lie algebras. Thus, the Lie algebraic methods here fail to
inform on representations of the group $H$. This failure arises
directly from the fact that the local coefficient system in homology
is non-trivial, that is, the action of $F_2$ on $H_*(F_3)$ is
non-trivial. For details, see Proposition
\ref{prop:poison}.

The above discussion suggests the following.

\begin{ques} \label{ques:faithful}
Let $\Gamma$ be a group that fits into a split exact sequence
\[
\begin{CD}
1 @>>> F @>>> \Gamma @>>> G @>>> 1
\end{CD}
\]
with $F$ a finitely generated free group and $G$ a
group that admits a finite dimensional faithful linear
representation. If $G$ acts trivially on the homology of $F$, does
$\Gamma$ then admit a finite dimensional faithful linear
representation? 
\end{ques}

\section{Centers}\label{sec:centers}

The purpose of this section is to identify the center of the
descending central series Lie algebra $\gr_*(\Gamma_n)$ in the case
where $\Gamma_n$ belongs to class of poly-free groups which satisfy
the conditions \eqref{eq:conditions}.  Recall that such groups are
given inductively by extensions
\[
\begin{CD}
1 @>{}>>  F_{\alpha(n)}@>{j}>> {\Gamma}_n @>{p}>> {\Gamma}_{n-1}
@>{}>>  1,
\end{CD}
\]
where $\Gamma_0$ is a free group; $F_{\alpha(n)}$ is a free group of rank at least $2$;
the map $p: {\Gamma}_n  \to \Gamma_{n-1}$ is a split epimorphism;
and the action of ${\Gamma}_{n-1}$ on $H_*(F_{{\alpha}(n)})$ is
trivial. Denote the center of a Lie algebra $\mathfrak g$ by
$\C(\mathfrak g)$, and write $F_{\alpha(0)} = \Gamma_0$.

\begin{prop}\label{prop:cyclic centers}
Let ${\Gamma}_n$ be a poly-free group satisfying conditions
\eqref{eq:conditions}. Then the center of $\gr_*(\Gamma_n)$ is
cyclic. Furthermore, if  $F_{{\alpha}(k)}$
is of rank at least two for every $k$, $0 \le k \le n$, 
then both the center of $\gr_*(\Gamma_n)$
and the center of ${\Gamma}_n$ are trivial.
\end{prop}

\begin{proof}
By \cite{falk-rand}, there is a short exact sequence of Lie algebras
\[
\begin{CD}
0 @>{}>> \gr_*(F_{\alpha(n)})@>{\gr_*(j)}>> \gr_*(\Gamma_n) @>{\gr_*(p)}>>
\gr_*(\Gamma_{n-1}) @>{}>>  0,
\end{CD}
\]
where $\gr_*(j)$ and $\gr_*(p)$ are the maps induced by $j$ and $p$, and the
first homology of $\Gamma_k$ is free abelian for each $k$. Since
$\Gamma_0$ is free of rank at least one, the center of
$\gr_*({\Gamma}_0)$ is either infinite cyclic if $\Gamma_0=\bbZ$, or is
trivial if the rank of $\Gamma_0$ is at least two \cite{copr}.

Assume inductively that the center of $\gr_*(\Gamma_{n-1})$ is
concentrated in degree $1$ and is cyclic.  Observe that
$\gr_*(p)(\C(\gr_*(\Gamma_n))$, the image of the center of
$\gr_*(\Gamma_n)$ in $\gr_*(\Gamma_{n-1})$, is also cyclic.  In
addition, it may be assumed inductively that the center
${\C}(\gr_*({\Gamma}_{n-1}))$ is trivial in case the rank of
$F_{\alpha(k)}$ is at least two for all $k$, $0 \le k \le n-1$.

Choose an element $\Delta_1$ in ${\C}(\gr_*({\Gamma}_{n}))$ which
projects to a generator of $\gr_*(p)(\C(\gr_*(\Gamma_n))$. Let $\Delta_2$
be an arbitrary element in the center ${\C}(\gr_*({\Gamma}_{n}))$.
If $\Delta_2$ is concentrated in degree greater than $1$, then
$\gr_*(p)(\Delta_2) = 0$. Consequently, there is an element $X$ in
$\gr_*(F_{\alpha(n)})$ such that $\gr_*(j)(X) = \Delta_2$. Since the
center of $\gr_*(F_{\alpha(n)})$ is trivial whenever the cardinality
of $\alpha(n)$ is at least two, $\Delta_2 = 0$. Thus it suffices to
assume that $\Delta_2$ is concentrated in degree $1$ 
in case the cardinality of $\alpha(n)$ is at least two.

It follows that there is an integer $q$ such that $q\Delta_1$ and
$\Delta_2$ project to the same element in
${\C}(\gr_*({\Gamma}_{n-1}))$. Then $\Delta = q\Delta_1 -\Delta_2$
is an element of ${\C}(\gr_*({\Gamma}_{n}))$, and is in the image of
$\gr_*(j)$, say $\Delta=\gr_*(j)(\delta)$, where $\delta \in
\C(\gr_*(F_{\alpha(n)}))$.  However, the center of
$\gr_*(F_{\alpha(n)})$ is trivial whenever the cardinality of
$\alpha(n)$ is at least two and so $q\Delta_1 -\Delta_2 = 0$. Thus
the center ${\C}(\gr_*({\Gamma}_{n}))$ is cyclic, and is trivial in
the case where the rank of $F_{\alpha(k)}$ is at least $2$ for all
$k$, $0 \le k \le n$.

To finish, assume that the rank of $F_{\alpha(k)}$ is at least $2$
for all $k$, $0 \le k \le n$. Since the center of $F_{\alpha(k)}$ is
trivial for $0 \le k \le n$, the center of ${\Gamma}_{n}$ is trivial
by the natural induction on $n$.
\end{proof}

The previous proposition suggests the following.

\begin{ques}
Let $\Gamma$ be a group which fits into an exact sequence
\[
\begin{CD}
1 @>>> F @>>> \Gamma @>>> G @>>> 1
\end{CD}
\]
where $F$ is a free group of rank at least two and the center of
$\gr_*(G)$ is cyclic. Is then the center of
$\gr_*(\Gamma)$ cyclic?  
Note that no assumptions are made regarding the existence of a splitting, 
or the action of $G$ on the homology of $F$.
\end{ques}

\section{Upper Triangular McCool Groups}

Let $F_n$ be the free group on the set $\{x_1, x_2, \cdots, x_n\}$.
Let $\IA_n$ denote the kernel of the natural map
\[
\Aut(F_n) \to \GL(n, \mathbb Z)
\]
induced by the map of $F_n$ to its abelianization.
As shown by Nielsen, and Magnus, the group $\IA_n$ is generated by
automorphisms $\beta_{i,j}$, $1\le i,j\le n$, $i\neq j$, and
$\Theta_{j;[s,t]}$, where $1\le j,s,t\le n$ and $j,s,t$ are
distinct, see \cite{magnus-karras-solit}.  These automorphisms are
given by
\[
\beta_{i,j}(x_k^{})=
\begin{cases}
x_k^{}, &\text{if $k \neq j$,}\\
x_i^{-1} x_j^{} x_i^{}, &\text{if $k=j$,}
\end{cases}
\quad\text{and}\quad \Theta_{j;[s,t]}(x_k^{})=
\begin{cases}
x_k^{}, &\text{if $k \neq j$,}\\
x_j^{} [x_s^{},x_t^{}], &\text{if $k=j$.}
\end{cases}
\]

Let $P\varSigma_n$ be the group of \emph{basis-conjugating
automorphisms}, the subgroup of $\IA_n$ generated by the
automorphisms $\beta_{i,j}$.  McCool \cite{MC} showed that the group
$P\varSigma_n$ admits a presentation with these generators and
relations
\begin{equation} \label{eq:McCool}
\left\{\begin{array}{ll}
{[}\beta_{i,j},\beta_{k,l}]& \text{for $i,j,k,l$ distinct}\\
{[}\beta_{i,j},\beta_{i,k}]& \text{for $i,j,k$ distinct}\\
{[}\beta_{j,k},\beta_{i,j}\beta_{i,k}]& \text{for $i,j,k$ distinct}\\
\end{array}\right\}.
\end{equation}

Let $P\varSigma_n^+$ be the \emph{upper triangular McCool group},
the subgroup of $P\varSigma_n$ generated by the automorphisms
$\beta_{i,j}$ with $i<j$.  The (relevant) relations \eqref{eq:McCool}
may be used to show that $P\varSigma_n^+$ is a semidirect product,
$P\varSigma_n^+ = F_{n-1} \rtimes P\varSigma_{n-1}^+$, where
$F_{n-1}$ is the free group on the set
\[
\{{\beta}_{i, n}: \; i = 1, 2, \dots n - 1\}.
\]
In other words, there is a split exact sequence
\[
\begin{CD}
1 @>>> F_{n-1} @>>> P\varSigma_n^+ @>>> P\varSigma_{n-1}^+ @>>> 1,
\end{CD}
\]
see \cite{copa} (where the automorphisms $\beta_{i,j}$ are denoted by ${\chi}(j,i)$).  
It is readily checked that the action of
$P\varSigma_{n-1}^+$ on $H_*(F_{n-1})$ is trivial.

The structure of the Lie algebra $\gr_*(P\varSigma_n^+)$ is
determined in \cite{copa}.  Denote the image of $\beta_{i,j}$ in
$\gr_*(P\varSigma_n^+)$ by $B_{i,j}$.  Then
$\gr_*(P\varSigma_n^+)$ has generators $B_{i,j}$, $1\le i<j\le
n$, and relations
\begin{subequations} \label{eq:McCool lie}
\begin{align}
{[}B_{i,j},B_{k,l}]&=0, \quad \text{for $i,j,k,l$ distinct,}
\label{first}\\
 [B_{i,j},B_{i,k}]&=0, \quad \text{for $i,j,k$ distinct,}
\label{second}\\
 [B_{j,k},B_{i,j}+B_{i,k}]&=0, \quad \text{for $i,j,k$ distinct.}
\label{third}
 \end{align}
\end{subequations}
Note that the last of these is one of the classical
infinitesimal braid relations.  It is not, however, the case that all of 
the infinitesimal braid relations are satisfied.  For example, if $i<j<k$, 
the element $[B_{i,k},B_{i,j}+B_{j,k}]$ need not vanish in 
$\gr_*(P\varSigma_n^+)$.

For $2\le k\le n$, let $V_k =\{B_{i,k}: \; 1\le i\le k - 1\}$,
and let $L[V_k]$ be the corresponding free Lie algebra.  There is an
additive isomorphism
\[
\gr_*(P\varSigma_n^+) \cong \bigoplus_{2\le k \le n} L[V_k].
\]
Moreover, it is clear from the relations \eqref{eq:McCool lie} that
the Lie algebra $\gr_*(F_{n-1})=L[V_n]$ is a Lie ideal of
$\gr_*(P\varSigma_n^+)$.

The next result is the portion of the Theorem \ref{thm-main}
pertaining to Case \textbf{I}, where $\Gamma_n = P\varSigma_n^+$.

\begin{thm}\label{prop-mccool}
The centralizer of $L[V_n]$ in ${\gr}_*(P\varSigma_n^+)$ is the
linear span of the element
\[{B}_n = \sum_{j=2}^n{B}_{1,j}.
\]
\end{thm}

\begin{proof}
Proposition \ref{prop:cyclic centers} implies that an element in the
centralizer of $L[V_n]$ must be of weight $1$. Let
\[
X = \sum_{1 \le i < j \le n}a_{i,j}{B}_{i,j}
\]
be an element that centralizes $L[V_n]$. Since $X$ centralizes
$L[V_n]$,
\[
[X, {B}_{1,n}] = 0.
\]
Relation \eqref{first} implies that
\[
\sum_{j=2}^{n-1} a_{1,j}[{B}_{1,j}, {B}_{1,n}] +
\sum_{i=2}^{n-1} a_{i,n}[{B}_{i,n}, {B}_{1,n}] = 0.
\]
Relation \eqref{second} implies that each term in the first sum
vanishes. Therefore,
\[
\sum_{i=2}^{n-1} a_{i,n}[{B}_{i,n}, {B}_{1,n}] = 0.
\]
But the last equation takes place in the free Lie algebra $L[V_n]$.
Thus, $a_{i,n} = 0$ for $1 < i < n$, and the element $X$ has the
form
\[
X = \sum_{1 \le i < j < n}a_{i,j}{B}_{i,j} +
a_{1,n}{B}_{1,n}.
\]

Fix $j$, $1 < j < n$. Start with the identity
\[
[X, {B}_{j,n}] = 0.
\]
Relation \eqref{first} implies that
\[
\sum_{i=1}^{j-1}a_{i,j}[{B}_{i,j}, {B}_{j,n}] +
\sum_{k=j+1}^{n-1} a_{j,k}[{B}_{j,k}, {B}_{j,n}] +
a_{1,n}[{B}_{1,n}, {B}_{j,n}] = 0.
\]
Relation \eqref{second} implies that each term in the middle sum
vanishes. Thus,
\[
\sum_{i=1}^{j-1}a_{i,j}[{B}_{i,j}, {B}_{j,n}] +
a_{1,n}[{B}_{1,n}, {B}_{j,n}] = 0.
\]
Using the relation \eqref{third}, notice that the next relation
follows
\[
-\sum_{i=1}^{j-1}a_{i,j}[{B}_{i,n}, {B}_{j,n}] +
a_{1,n}[{B}_{1,n}, {B}_{j,n}] = 0
\]
which implies that
\[
-\sum_{i=2}^{j-1}a_{i,j}[{B}_{i,n}, {B}_{j,n}] + (-a_{1,j} +
a_{1,n})[{B}_{1,n}, {B}_{j,n}] = 0.
\]
Since all the commutators are linearly independent, $a_{i,j} = 0$
for $1 < i < j$, and $a_{1,j} = a_{1,n}$. Since $j$ was arbitrary,
these relations hold for all $j$ with $1 < j < n$. The result
follows.
\end{proof}

\section{Surface Group Orbit Configuration Spaces}
Let $M$ be a manifold without boundary, and let $G$ be a discrete
group which acts freely on $M$. The {\em orbit configuration space}
consists of all ordered $n$-tuples of points in $M$ which lie in
distinct orbits:
\[
\conf^G(M,n)=\{(x_1,\dots,x_n)\in M^n:\; G\cdot x_i\cap G \cdot x_j=
\emptyset\ \text{if}\ i\neq j\}.
\]

In \cite{xi}, Xicot\'encatl proves that, for $\ell \le n$,
projection onto the first $\ell$-coordinates,
\[
p_G:\conf^G(M,n) \to \conf^G(M,\ell),
\]
is a locally trivial bundle, with fiber $\conf^G(M\setminus
Q_\ell^G,n-\ell)$, where $Q_{\ell}^G$ denote the union of $\ell$
distinct orbits, $G\cdot x_1,\dots,G\cdot x_n$, in $M$. This result 
generalizes the Fadell-Neuwirth theorem \cite{dane}.

Let $G$ be a discrete subgroup of $PSL(2,{\bbR})$ acting freely and
properly discontinuously on the upper--half plane ${\bbH}^2$ by
fractional linear transformations. Let ${\conf}^{G}({\bbH}^2, n)$ be
the corresponding orbit configuration space and $P_G(n)$ its
fundamental group. Then there is a fiber bundle
\[
{\bbH}^2 {\setminus} Q_{n-1}^G \to \text{\conf}^G({\bbH}^2, n) \to
\conf^G({\bbH}^2, n - 1).
\]
This induces a split exact sequence:
\[
\begin{CD}
1 @>>> F @>>> P_G(n) @>>> P_G(n - 1) @>>> 1,
\end{CD}
\]
where $F$ is the free group on $Q_{n-1}^G$. Notice that in this case
$F$ is an infinitely generated free group.

The structure of ${\gr}_*(P_G(n))$ is determined in \cite{ckx}. The
generators are $B_{i,j}^{\sigma}$, where $1 \le i < j \le n$ and
$\sigma \in G$. In this case, the infinitesimal braid relations are:
\begin{equation} \label{eq:inf surface}
\begin{aligned}
{[}B_{i,j}^{\sigma}, B_{s,t}^{\tau}] &= 0, \quad \text{if $\{i,
j\}{\cap}\{s, t\}
= \emptyset$,}\\
[B_{i,j}^{\tau}, B_{s,j}^{{\tau}{\sigma}^{-1}} + B_{i,s}^{\sigma}]
&= 0, \quad \text{if $1 \le i < s < j \le n$,}\\
[B_{i,s}^{\sigma}, B_{i,j}^{\tau} + B_{s,j}^{{\tau}{\sigma}^{-1}}]
&= 0, \quad \text{if $1 \le i < s < j \le n$,}\\
[B_{s,j}^{{\tau}{\sigma}^{-1}}, B_{i,j}^{\tau} + B_{i,s}^{\sigma}]
&= 0, \quad \text{if $1 \le i < s < j \le n$.}
\end{aligned}
\end{equation}

The last equation follows from the previous two equations. Note
that, for $0 < i < k < n$, this last equation implies that
\begin{equation} \label{eq:number 1}
[B_{k,n}^{{\tau}{\sigma}^{-1}}, B_{i,n}^{\tau} + B_{i, k}^{\sigma}]
= 0.
\end{equation}

Let $V_j$ be the linear span of the set $\{B_{i,j}^{\sigma}:\; 1 \le
i < j,\; {\sigma}\in G\}$. Then, there is a splitting
\[
{\gr}_i(P_G(n)) = \bigoplus_{j=2}^n{\gr}_i(L[V_j]),
\]
as abelian groups. 

The next result is the portion of the Theorem \ref{thm-main}
pertaining to Case \textbf{II}, where $\Gamma_n = P_G(n)$.

\begin{thm}\label{thm-orbit}
The centralizer of $L[V_n]$ in ${\gr}_*(P_G(n))$ is zero:
\[
{\C}_{{\gr}_*(P_G(n))}(L[V_n]) = 0.
\]
\end{thm}

\begin{proof}
Proposition \ref{prop:cyclic centers} implies that elements in the
centralizer are of weight  $1$. Let
\[x = \sum_{k = 1}^{\ell}a_kB_{i_k,j_k}^{{\sigma}_k} +
\sum_{m=1}^pb_mB_{m,n}^{{\tau}_m}
\]
be such an element.  Since $G$ is infinite, there exists $\tau \in
G$ such that
\[
{\tau}{\sigma}_k^{-1} \not= {\tau}_m, \quad \text{for all}\; \; k =
1, 2, \dots , {\ell}, \; m = 1, 2, \dots , p.
\]
Since $x$ centralizes $L[V_n]$,
\[
0 = [x, B_{1,n}^{\tau}] = \sum_{k =
1}^{\ell}a_k[B_{i_k,j_k}^{{\sigma}_k}, B_{1,n}^{\tau}] +
\sum_{m=1}^pb_m[B_{m,N}^{{\tau}_m}, B_{1,n}^{\tau}].
\]
The infinitesimal braid relations \eqref{eq:inf surface} imply that
\[
[B_{i_k,j_k}^{{\sigma}_k}, B_{1,n}^{\tau}] = 0, \quad \text{if}\;\;
1 \not=i_k,
\]
and that
\[
[B_{1,j_k}^{{\sigma}_k}, B_{1,n}^{\tau}] =
-[B_{j_k,n}^{{\tau}{\sigma}_k^{-1}}, B_{j_k,n}^{\tau}].
\]
Thus, for eack $1 \le i_k < n$,
\[
[B_{i_k,j_k}^{{\sigma}_k}, B_{1,n}^{\tau}] =
-[B_{j_k,n}^{{\tau}{\sigma}_k^{-1}}, B_{1,n}^{\tau}],
\]
and
\[
-\sum_{k = 1}^{\ell}a_k[B_{j_k,n}^{{\tau}{\sigma}_k^{-1}},
B_{1,n}^{\tau}] + \sum_{m=1}^pb_m[B_{m,n}^{{\tau}_m},
B_{1,n}^{\tau}] = 0.
\]
Since the commutators are all different and they are linearly
independent, all the coefficients are equal to $0$ and thus $x = 0$.
\end{proof}

\section{Cyclic Group Orbit Configuration Spaces}\label{sec:cyclic orbit}
Let $G={\bbZ}/r{\bbZ}$ be a finite cyclic group.  The group $G$ acts freely 
on the manifold $M={\bbC}^* =
{\bbC}\setminus\{0\}$ by multiplication by the primitive $r$-th root
of unity $\zeta=\exp(2\pi\ii/r)$. The corresponding orbit configuration space is given by
\[
\conf^G({\bbC}^*,n)=\{(x_1,\dots,x_n) \in ({\bbC}^*)^n \mid x_j \neq
\zeta^p x_i \ \text{for $i\neq j$ and $1\le p \le r$}\}.
\]
Denote the fundamental group of $\conf^G({\bbC}^*,n)$ by $P(r,n)$.  

Let $\ell_n = rn+1$, and define a map $g_n: \conf^G(\bbC^*,n) \to
\conf(\bbC,\ell_n)$ from the orbit configuration space to the
classical configuration space by sending a point to its orbits
(together with $0$).  Explicitly, if $(x_1,\dots,x_n) \in
\conf^G(\bbC^*,n)$, define
\[
g_n(x_1,\dots,x_n) = (0,x_1,\zeta x_1,\dots,\zeta^{r-1}
x_1,x_2,\zeta x_2,\dots,\zeta^{r-1} x_2, \dots\dots,x_n,\zeta
x_n,\dots,\zeta^{r-1} x_n)
\]
in $\conf(\bbC,\ell_n)$. Then, one has the following result (see
\cite[Thm.~2.1.3]{mono} and \cite[\S{3}]{ccx}).

\begin{thm} \label{thm:orbit pullback}
The orbit configuration space bundle $p_G: \conf^G(\bbC^*,n+1) \to
\conf^G(\bbC^*,n)$ is equivalent to the pullback of the
configuration space bundle $p: \conf(\bbC,\ell_n+1) \to
\conf(\bbC,\ell_n)$ along the map $g_n$.
\end{thm}

Passing to fundamental groups, there is an induced commutative
diagram with split rows,
\begin{equation} \label{eqn:groups}
\begin{CD}
1 @>>> F_{\ell_n}    @>>> P(r,n+1)   @>>> P(r,n) @>>>1\\
@.     @VV{\rm{id}}V      @VVV         @VV{(g_n)_{*}}V\\
1 @>>> F_{\ell_n}    @>>> P_{\ell_n+1} @>>> P_{\ell_n} @>>>1
\end{CD}
\end{equation}
where $F_N$ is a free group on $N$ generators.  Passing further to
descending central series Lie algebras, there is a commutative
diagram, again with split rows (see \cite[\S{4}]{ccx}).
\begin{equation} \label{eqn:algebras}
\begin{CD}
0 @>>> L[V_{\ell_n}]    @>>> {\gr}_*(P(r,n+1))   @>>>
{\gr}_*(P(r,n))
@>>>0\\
@.     @VV{\rm{id}}V      @VVV         @VV{\gr_*(g_n)}V\\
0 @>>> L[V_{\ell_n}]    @>>> {\gr}_*(P_{\ell_n+1}) @>>>
{\gr}_*(P_{\ell_n}) @>>>0
\end{CD}
\end{equation}
This realizes the Lie algebra ${\gr}_*(P(r,n+1))$ as the semidirect
product of ${\gr}_*(P(r,n))$ by $L[V_{\ell_n}]$ determined by the
homomorphism
\[
\theta_{\ell_n} \circ \gr_*(g_n): {\gr}_*(P(r,n)) \to
{\rm{Der}}(L[V_{\ell_n}]),
\] 
where $\theta_N: {\gr}_*(P_N) \to
{\rm{Der}}(L[V_N])$ is given by
$\theta_N(B_{i,j})={\rm{ad}}(B_{i,j})$, see \cite[Thm.~4.4]{ccx}.
More explicitly, the structure of ${\gr}_*(P(r,n))$ is given in the
following theorem.

\begin{thm}[\cite{mono}] \label{thm:PureMonomialBrackets}
Let ${\gr}_*(P(r,n))$ be the Lie algebra associated to the descending 
central series of the group $P(r,n)=\pi_1(\conf^G(\bbC^*,n)$, where 
$G=\bbZ/r\bbZ$.  Then,
\[
{\gr}_*(P(r,n)) \cong \bigoplus_{j=0}^{n-1} L(rj+1)
\]
as abelian groups, where $L(rj+1)$ is generated by $Z_{j+1}$ and
$B_{i,j+1}^{(p)}$, $1\le i\le j$, $1\le p \le r$.  The Lie bracket
relations in ${\gr}_*(P(r,n))$ are given by
\[
\begin{aligned}
\bigl[Z_j+Z_l+B_{j,l}^{(1)}+B_{j,l}^{(2)}+\dots+B_{j,l}^{(r)},Y\bigr]
&=0 \quad \text{for $Y=Z_l$, $Y=B_{j,l}^{(p)}$, $1\le p\le r$,} \\
[B_{i,j}^{(p)}+B_{i,k}^{(q)}+B_{j,k}^{(m)},Y]&=0 \quad\text{for
$Y=B_{i,k}^{(q)},B_{j,k}^{(m)}$, $q\equiv p+m \mod
r$,}\\
[B_{i,j}^{(p)},B_{k,l}^{(q)}]&=0 \quad\text{for
$\{i,j\}\cap\{k,l\}=\emptyset$, $1\le p,q\le r$,
and}\\
[Z_{k},B_{i,j}^{(p)}]&=0 \quad\text{for $k\neq i,j$ and $1\le p \le
r$.}
\end{aligned}
\]
\end{thm}

\begin{prop} \label{prop:injections}
The map $g_n:\conf^G(\bbC^*,n) \to \conf(\bbC,\ell_n)$ induces
injections on fundamental groups and descending central series Lie
algebras. More precisely, the maps
\[
(g_n)_{*}: P(r,n) \longrightarrow P_{{\ell}_n} \quad\text{and}\quad \gr_*(g_n): {\gr}_*(P(r,n))
\longrightarrow {\gr}_*(P_{\ell_n})
\]
induced by $g_n$, are monomorphisms.
\end{prop}

\begin{proof}
The proof is by induction on $n$.  In the case $n=1$, notice that
$\conf^G(\bbC^*,1)=\bbC^*$, and for $x\in\bbC^*$, $g_1(x)=(0,x,\zeta
x,\dots,\zeta^{r-1}x) \in \conf(\bbC,r+1)$.  Let $\gamma \in
\pi_1(\bbC^*)$ be the standard generator, and check that
$(g_1)_*(\gamma)=\Delta(r+1)$ generates the center of $P_{r+1}$. It
follows that both $(g_1)_*:P(r,1) \to P_{\ell_1}$ and $\gr_*(g_1):
{\gr}_*(P(r,1)) \to {\gr}_*(P_{\ell_1})$ are injective.

Assume inductively that $(g_n)_*:P(r,n) \to P_{\ell_n}$ and
$\gr_*(g_n): {\gr}_*(P(r,n)) \to {\gr}_*(P_{\ell_n})$ are injective. It
must be shown that $(g_{n+1})_*:P(r,n+1) \to P_{\ell_{n+1}}$ and
$\gr_*(g_{n+1}): {\gr}_*(P(r,n+1)) \to {\gr}_*(P_{\ell_{n+1}})$ are
also injective, where $\ell_n=rn+1$ and $\ell_{n+1}=r(n+1)+1$. Let
$\tilde{g}_n:\conf^G(\bbC^*,n+1) \to \conf(\bbC,\ell_n+1)$ denote
the map on the pullback induced by $g_n$. Note that
\[
\tilde{g}_n(x_1,\dots,x_n,z)= (0,x_1,\zeta x_1,\dots,\zeta^{r-1}
x_1, \dots\dots,x_n,\zeta x_n,\dots,\zeta^{r-1} x_n, z).
\]
The map $\tilde{g}_n$ may be factored as follows.  Let
$p_{m,k}:\conf(\bbC,m) \to \conf(\bbC,k)$ be the projection which
forgets the last $m-k$ points.  Then $\tilde{g}_n = p_{m,k} \circ
g_{n+1}$, where $m=\ell_{n+1}$ and $k=\ell_n + 1$.

Since $(g_n)_*$ and $\gr_*(g_n)$ are injective by induction, it
follows from \eqref{eqn:groups} and \eqref{eqn:algebras} that
$(\tilde{g}_n)_*$ and $\gr_*(\tilde{g}_n)$ are also injective.  This,
together with the fact that $\tilde{g}_n = p_{m,k} \circ g_{n+1}$,
implies that $(g_{n+1})_*$ and $\gr_*(g_{n+1})$ are injective.
\end{proof}

\begin{cor} \label{cor:cyclic orbit linear}
The group $P(r,n)$ is linear.
\end{cor}
\begin{proof}
The group $P(r,n)$ embeds in the pure braid group, which is linear
\cite{bigelow, krammer}.
\end{proof}

The next result is the portion of the Theorem \ref{thm-main}
pertaining to Case \textbf{III}.  It is notationally convenient to
state the result for the group $\Gamma_{n+1}=P(r,n+1)$.

\begin{thm}
The centralizer of $L[V_{\ell_n}]$ in ${\gr}_*(P(r,n+1))$ is the
linear span of the element
\[
\Delta(r,n+1) = \sum_{k=1}^{n+1} Z_k +\sum_{p=1}^r \sum_{1\le i<j
\le n+1} B_{i,j}^{(p)}.
\]
\end{thm}

\begin{proof}
Denote the generators of $L[V_{\ell_n}]$ by $Z_{n+1}$ and
$B_{i,n+1}^{(p)}$, $1\le i \le n$, $1\le p \le r$, and let
\[
\B_{\ell_n}= Z_{n+1}+\sum_{p=1}^r \sum_{i=1}^n B_{i,n+1}^{(p)}.
\]
Let $x \in {\gr}_*(P(r,n+1))$, and assume that $[x,B]=0$ for all
$B\in L[V_{\ell_n}]$.  Write $x=u+v$, where $u \in {\gr}_*(P(r,n))$
and $v \in L[V_{\ell_n}]$.  Then, for all $B \in L[V_{\ell_n}]$, it
follows that $\gr_*(\tilde{g}_n)[x,B]=[\gr_*(g_n)(u)+v,B]=0$ in
${\gr}_*(P_{\ell_n+1})$. So $\gr_*(g_n)(u)+v = k\cdot \Delta(\ell_n+1)$
for some constant $k$.  Consequently, $\gr_*(g_n)(u) = k \cdot
\Delta(\ell_n)$ and $v=k\cdot \B_{\ell_n}$.  Since $\gr_*(g_n):
{\gr}_*(P(r,n)) \to {\gr}_*(P_{\ell_n})$ is injective (as is
$\gr_*(\tilde{g}_n)$), it follows that the centralizer of
$L[V_{\ell_n}]$ in ${\gr}_*(P(r,n+1))$ is the linear span of the
element
\[
\gr_*(\tilde{g}_n)^{-1}(\Delta(\ell_n+1)) =
\gr_*(g_n)^{-1}(\Delta(\ell_n))+ \B_{\ell_n}.
\]

So it suffices to show that
\[
\gr_*(g_n)^{-1}(\Delta(\ell_n))= \sum_{k=1}^{n} Z_k +\sum_{p=1}^r
\sum_{1\le i<j \le n} B_{i,j}^{(p)}.
\]
The map $g_n:\conf^G(\bbC^*,n) \to \conf(\bbC,\ell_n)$ is the
restriction of the affine transformation $g_n:\bbC^n \to
\bbC^{\ell_n}$, defined by the same formula, and, abusing notation,
denoted by the same symbol. The orbit configuration space
$\conf^G(\bbC^*,n)$ may be realized as the complement of the
hyperplane arrangement $\A$ in $\bbC^n$ with hyperplanes $H_i =
\{x_i=0\}$, $1\le i \le n$, and $H_{i,j}^{(p)}=\{x_i = \zeta^p
x_j\}$, $1\le i < j \le n$, $1\le p \le r$.  The generators of the
Lie algebra ${\gr}_*(P(r,n))$ are in one-to-one correspondence with
the hyperplanes of $\A$.  If $B_H \in {\gr}_*(P(r,n))$ denotes the
generator corresponding to $H \in \A$, it follows from
\cite[Prop.~3.4]{ccx} that
\[
\gr_*(g_n)(B_H) = \sum_{g_n(H) \subset H_{r,s}} B_{r,s},
\]
where $H_{r,s} = \{x_r=x_s\} \subset \bbC^{\ell_n}$, and $B_{r,s}
\in {\gr}_*(P_{\ell_n})$ is the corresponding generator of the
descending central series Lie algebra of the pure braid group.  Let
$S_i = \{ H_{r,s} \mid g_n(H_i) \subset H_{r,s}\}$ and
$S_{i,j}^{(p)}= \{H_{r,s} \mid g_n(H_{i,j}^{(p)}) \subset
H_{r,s}\}$.  Then, one can check that the sets $S_i$, $1\le i \le
n$, and $S_{i,j}^{(p)}$, $1\le i < j \le n$, $1\le p \le r$, form a
partition of the (entire) set of $\binom{\ell_n}{2}$ hyperplanes
$H_{r,s}$ in $\bbC^{\ell_n}$.  It follows that
\[
\gr_*(g_n)^{-1}(\Delta(\ell_n))= \gr_*(g_n)^{-1}\Big(\sum_{1\le r < s \le
\ell_n} B_{r,s}\Big)= \sum_{k=1}^{n} Z_k +\sum_{p=1}^r \sum_{1\le
i<j \le n} B_{i,j}^{(p)},
\]
completing the proof.
\end{proof}

\section{Fiber-Type Arrangements}\label{sec:fiber-type}

In light of Corollary \ref{cor:cyclic orbit linear}, it is natural
to speculate that the fundamental group of the complement of an
arbitrary fiber-type hyperplane arrangement is linear.  The purpose
of this section is to show that this is indeed the case.

A hyperplane arrangement $\A$ is a finite collection of codimension
one affine subspaces of Euclidean space $\bbC^{n}$. See Orlik and
Terao \cite{OT} as a general reference on arrangements. The
complement of an arrangement $\A$ is the manifold
$X=X(\A)=\bbC^{n}\setminus\bigcup_{H\in\A}H$.  Denote the
fundamental group of the complement by $G(\A)=\pi_1(X(\A))$.

\begin{defn} \label{defn:slfdef}
A hyperplane arrangement $\A$ in $\bbC^{n+1}$ is \emph{strictly
linearly fibered} if there is a choice of coordinates
$(\b{x},z)=(x_1,\dots,x_n,z)$ on $\bbC^{n+1}$ so that the
restriction, $p$, of the projection $\bbC^{n+1}\to\bbC^n$,
$(\b{x},z)\mapsto \b{x}$, to the complement $X(\A)$ is a fiber
bundle projection, with base $p(X(\A))=X(\B)$, the complement of an
arrangement $\B$ in $\bbC^n$, and fiber the complement of finitely
many points in $\bbC$. In this case $\A$ is said to be strictly
linearly fibered over $\B$.
\end{defn}

Let $\A$ be an arrangement in $\bbC^{n+1}$, strictly linearly
fibered over $\B \subset \bbC^n$. For each hyperplane $H$ of $\A$,
let $f_H$ be a linear polynomial with $H=\ker f_H$.  Then
$Q(\A)=\prod_{H\in\A} f_H$ is a defining polynomial for $\A$.  From
the definition, there is a choice of coordinates for which a
defining polynomial for $\A$ factors as $Q(\A) = Q(\B)\cdot
\phi(\b{x},z)$, where $Q(\B)=Q(\B)(\b{x})$ is a defining polynomial
for $\B$, and $\phi(\b{x},z)$ is a product
\[
\phi(\b{x},z) = (z-g_1(\b{x}))(z-g_2(\b{x}))\cdots (z-g_m(\b{x})),
\]
with $g_j(\b{x})$ linear.  Define $g:\bbC^{n}\to\bbC^{m}$ by
\begin{equation} \label{eq:rootmap}
g(\b{x})=\bigl(g_1(\b{x}),g_2(\b{x}),\dots,g_m(\b{x})\bigr),
\end{equation}
Since $\phi(\b{x},z)$ necessarily has distinct roots for any
$\b{x}\in X(\B)$, the restriction of $g$ to $X(\B)$ takes values in
the configuration space $\conf(\bbC, m)$. The next result
generalizes Theorem \ref{thm:orbit pullback}.

\begin{thm}[\cite{mono}]
\label{thm:slfpullback} Let $\B$ be an arrangement of $m$
hyperplanes, and let $\A$ be an arrangement of $m+n$ hyperplanes
which is strictly linearly fibered over $\B$.  Then the bundle
$p:X(\A)\to X(\B)$ is equivalent to the pullback of the bundle of
configuration spaces $p_{n+1}:\conf(\bbC,n+1)\to \conf(\bbC,n)$
along the map $g$. Consequently, the bundle $p:X(\A)\to X(\B)$
admits a cross-section, and has trivial local coefficients in
homology.
\end{thm}

Passing to fundamental groups, there is an induced commutative
diagram with split rows.
\begin{equation} \label{eqn:slf groups}
\begin{CD}
1 @>>> F_{n}    @>>> G(\A)   @>>> G(\B) @>>>1\\
@.     @VV{\rm{id}}V      @VVV         @VV{g_{*}}V\\
1 @>>> F_{n}    @>>> P_{n+1} @>>> P_{n} @>>>1
\end{CD}
\end{equation}
realizing $G(\A)$ as a pullback.

\begin{lem} \label{lem:linear pullback}
Given a pullback of groups
\[
\begin{CD}
@. @. P   @>>> H \\
@.     @.      @VVV         @VVV\\
1 @>>> K    @>>> G @>>> Q @>>>1
\end{CD}
\]
if $G$ and $H$ are linear, then the pullback $P$ is also linear.
\end{lem}

\begin{proof}
The pullback $P$ is a subgroup of $G \times H$, which is linear if
$G$ and $H$ are.
\end{proof}

\begin{cor} \label{cor:slf linear}
If $\A$ is strictly linearly fibered over $\B$, and $G(\B)$ is
linear, then $G(\A)$ is also linear.
\end{cor}

\begin{defn} \label{defn:ftdef}
An arrangement $\A=\A_1$ of finitely many points in $\bbC^1$ is
\emph{fiber-type}.  An arrangement $\A=\A_n$ of hyperplanes in
$\bbC^n$ is \emph{fiber-type} if $\A$ is strictly linearly fibered
over a fiber-type arrangement $\A_{n-1}$ in $\bbC^{n-1}$.
\end{defn}

Examples include the braid arrangement with defining polynomial
$Q(\A)=\prod_{i<j}(y_i-y_j)$, and complement $X(\A)=\conf(\bbC,n)$,
and the full monomial arrangement with defining polynomial
$Q(\A)=x_1\cdots x_n \prod_{i<j}(x^r_i-x^r_j)$, and complement
$X(\A)=\conf^G(\bbC^*,n)$, where $G=\bbZ/r\bbZ$.

\begin{thm} \label{thm:ft linear}
The fundamental group of the complement of a fiber-type hyperplane
arrangement is linear.
\end{thm}
\begin{proof}
Let $\A=\A_n$ be a fiber-type arrangement in $\bbC^n$, with
complement $X(\A_n)$.  The proof is by induction on $n$.

In the case $n=1$, denote the cardinality of $\A_1$ by $d$.  If
$d=0$ ($\A_1$ is the empty arrangement), then $X(\A_1)=\bbC$ and
$G(\A_1)$ is the trivial group.  If $d>0$, then $X(\A_1)$ has the
homotopy type of a bouquet of $d$ circles and $G(\A_1)=F_d$ is a
free group on $d$ generators.  So $G(\A_1)$ is linear.

Assume the result holds for any fiber-type arrangement $\A_n$ in
$\bbC^n$, and let $\A_{n+1}$ be a fiber-type arrangement in
$\bbC^{n+1}$.  Then $\A=\A_{n+1}$ is strictly linearly fibered over
$\B=\A_n$, a fiber-type arrangement in $\bbC^n$.  By induction, the
fundamental group $G(\B)=\pi_1(X(\B))$ is linear.  Hence,
$G(\A)=\pi_1(X(\A))$ is also linear by Corollary \ref{cor:slf
linear}.
\end{proof}

\section{The Poison Group} \label{sec:poison}

One group which does not admit a faithful finite dimensional linear
representation is the so-called poison group.  This group appears in the work 
of Formanek and Procesi \cite{fopr}, where it is realized as a subgroup of 
$\Aut(F_n)$ for $n \ge 3$, proving that the latter admits no faithful finite 
dimensional linear representation.  The poison group may also be realized as a 
subgroup of $\IA_n$ for $n \ge 5$, see Pettet \cite{pettet}.  On the other hand, 
Brendle and Hamidi-Tehrani \cite{BHT} have shown that the mapping class group
of a genus $g$ surface with one fixed point has no subgroup isomorphic to the poison group.

The purpose of this section to show
how the Lie algebraic criteria in this article fail in a strong way
for the poison group $H$, a group given by a split extension
\[
\begin{CD}
1 @>>>  F_3 @>{j}>> H @>{p}>> F_2 @>>> 1.
\end{CD}
\]
The group $H$ admits the following presentation:
\begin{equation} \label{eq:poison pres}
H=\langle a_1,a_2,a_3,\phi_1,\phi_2 \mid \phi_i^{} a_j^{}
\phi_i^{-1}=a_j^{}, \phi_i^{} a_3^{} \phi_i^{-1}=a_3^{} a_i^{},
i,j=1,2\rangle.
\end{equation}
 From this presentation, it is clear that $H = F_3 \rtimes F_2$ is a
semidirect product, where
 $F_3$ is generated by $\{a_1,a_2,a_3\}$ and $F_2$ by $\{\phi_1,\phi_2\}$.
Thus, $H$ is poly-free, and it is natural to consider how the
structure of the descending central series Lie algebra fails to
inform on representations for this group. For a group $G$, let
$\alpha:G \to H_1(G)$ denote the abelianization map.

\begin{prop}\label{prop:poison}
There is a split exact sequence of Lie algebras
\[
\begin{CD}
0 @>>> \mathbb Z @>{}>>  {\gr}_*(H) @>{\gr_*(p)}>> {\gr}_*(F_2) @>>> 0
\end{CD}
\]
with the center, ${\C}({\gr}_*(H))$, given by $\mathbb Z$, 
generated the class of $a_3$.

The induced map
\[
\begin{CD}
  {\gr}_*(F_3) @>{\gr_*(j)}>>  {\gr}_*(H)
\end{CD}
\] factors through the center ${\C}({\gr}_*(H)) = \mathbb Z$, and
\[
\begin{CD}
  {\gr}_*(H)/{\C}({\gr}_*(H)) @>{\gr_*(j)}>> {\gr}_*(F_2)
\end{CD}
\] is an isomorphism of Lie algebras.

Furthermore, the natural map
\[p \times \alpha: H \longrightarrow F_2\times H_1(H)
\]
has non-trivial kernel, but induces a monomorphism
\[
\gr_*(p \times \alpha): {\gr}_*(H) \longrightarrow {\gr}_*(F_2\times H_1(H))
\]
on the level of descending central series Lie algebras.
\end{prop}

Consequently, the Lie algebra obtained from the descending central
series of $H$ provides little information about embeddings as the
subgroup $F_3=\langle a_1,a_2,a_3\rangle$ has image which factors
through $\mathbb Z$ on the level of Lie algebras.

\begin{proof}
It follows from the presentation \eqref{eq:poison pres} that, for
$1\le i,j\le 2$, the relations $[\phi_i,a_j]=1$ and
$[\phi_i,a_3]=a_3^{}a_i^{}a_3^{-1}$ hold in $H$.  Denote the images
of the generators $\phi_i$ and $a_j$ in $\gr_*(H)$ by the same
symbols.  Then, for $1\le i\le 2$, $a_i=0$ in $\gr_1(H)$ since $a_i$
is conjugate to a commutator in $H$.  For any element
$X \in \gr_*(H)$, it follows that $[a_i,X]=0$ in $\gr_*(H)$ if 
$1\le i\le 2$. Also, since $[\phi_i,a_3]=[a_3,a_i]\cdot a_i$ in $H$,
$[\phi_i,a_3]=0$ in $\gr_1(H)$. It follows that the map
$\gr_*(j):\gr_*(F_3) \to \gr_*(H)$ factors through the Lie subalgebra
generated by $a_3$:
\[
\begin{CD}
\gr_*(j): {\gr}_*(F[a_1,a_2,a_3]) @>>> L[a_3]
@>>> {\gr}_*(H)
\end{CD}
\]

Since $[\phi_i, a_3]$ is zero in ${\gr}_*(H)$, the class of $a_3$
centralizes ${\gr}_*(H)$. It follows that the Lie algebra kernel of
$\gr_*(p):\gr_*(H) \to \gr_*(F_2)$ is exactly $\C({\gr}_*(H))$, a copy of
the integers generated by $a_3$.

Note that $H_1(H)$ is free abelian of rank $3$, generated by the
classes of $\phi_1$, $\phi_2$, and $a_3$.  It follows that the
natural map $p \times \alpha: H \to F_2 \times H_1(H)$ induces a
monomorphism
\[
\gr_*(p \times \alpha): {\gr}_*(H) \longrightarrow {\gr}_*(F_2 \times H_1(H))
\]
on the level of Lie algebras, which completes the proof.
\end{proof}

\end{document}